\documentclass[12pt,a4paper]{article}
\usepackage[latin2]{inputenc}
\usepackage{graphicx}
\usepackage{amsmath}
\usepackage{hyperref}
\usepackage{amssymb}                %Use this in pdflatex mode  %obligatory package

\newtheorem{theorem}{Theorem}[section]

\newtheorem{remark}[theorem]{Remark}

\begin{document}

\begin{center}

%%\title{

{\Large Stability of a retrovirus dymanic model}
%{\Large Mathematical model of anticancer virotherapy with retrovirus
%}
%}
\bigskip

%{\large 
%A.K. and A.R. }

%\author{
{\large 
Andrei Korobeinikov$^1$ and Alexander  Rezounenko$^{2,3,}$ }\footnote[4]{Corresponding author: A. Rezounenko (email: rezounenko@gmail.com)}

\medskip

$^1$Centre de Recerca Matematica
Campus de Bellaterra, \\ 
Edifici C, 
08193 Bellaterra (Barcelona), Spain \\
$^{2}$ V.N.Karazin Kharkiv National University, 
  4 Svobody sqr.,\\  Kharkiv, 61022,  Ukraine\\
$^{3}$The Czech Academy of Sciences,\\  Institute of Information Theory and Automation,\\
 P.O. Box 18, 182\,08 Praha, CR

\end{center}
%\maketitle
%%%\address{Kharkiv}
\begin{abstract}
%%%%%%%%%%%%%%%%%%%\marginpar{!!!}%%%%%%%%%%%%%%%%%%%%%%%%%%%%%%%%%
A retrovirus dynamic model is proposed. We pay attention to the case when viral pathogenicity is low and the infected cells are able to reproduce. Using Lyapunov function method we study stability properties of an inner equilibrium of the model.  The equilibrium represents a chronic disease steady state. 

\end{abstract}

\marginpar{\tiny 2018.12.29.AR}  

2010 Mathematics Subject Classification: 
34D20; %%%Stability
93D05; %%%Lyapunov and other classical stabilities (Lagrange, Poisson, $L^p, l^p$, etc.)
92B05. %%General biology and biomathematics

\medskip

 Keywords: Evolution equations; Lyapunov stability; %state-dependent delay, 
 virus infection model;
 
  anticancer virotherapy.

\section{Introduction}\label{AR-2018-12-introduction}

 There are two major modes of viral replication: lytic and non-lytic. The dynamics of lytic virus is considered in many papers (see, e.g. \cite{Nowak-Bangham-S-1996,Perelson-Neumann-Markowitz-Leonard-Ho-S-1996}). Non-lutic is almost out of consideration.
In the nonlytic reproduction virus stays dormant in an infected cell, that reproduce. Virus stays dormant until the cell stars to exhibit signs of exhaustion. When the fist such signs appear, virus starts fast replication in the cell, kills it and brakes the cell membrane releasing virus particles.

The idea of virotherapy consist in  
using viruses in delivery  their genetic material (a piece of RNA or DNA) into host cells that need to be treated. There are several branches of virotherapy such as viral immunotherapy, anti-cancer oncolytic viruses and  viral vectors for gene therapy. %% \cite{}. 
%%%}%%%
Anticancer virotherapy is a new and promising method of anticancer therapy.

However, the majority of viruses that are used for virotherapy are of highly virulent and highly pathogenic type, and, therefore, the infected cells usually do not reproduce (producing the virus instead). Typical mathematical models reflect this fact.

The following classical model of virus dynamics was proposed in \cite{Nowak-Bangham-S-1996,Perelson-Neumann-Markowitz-Leonard-Ho-S-1996} 
\begin{equation}\label{AR-2018-12-01}
\left\{
\begin{array}{l}
\frac{dC(t)}{dt}=\lambda -dC(t)-\beta C(t)V(t), \\
\frac{d I(t)}{dt}=\beta C(t)V(t)-aI(t), \\
\frac{dV(t)}{dt}=aN\, I(t)-kV(t), \\
\end{array}
\right.
\end{equation}
where $C(t), I(t), V(t)$ represent the concentration (or the total number) of
non-infected host cells, infected cells and free virions at time $t$, respectively. All the constants $\lambda, d, \beta, a, N, k$ in (\ref{AR-2018-12-01}) are positive. 
The non-infected cells are produced at rate $\lambda$, die at rate $d$ and become infected at rate $\beta$. Infected cells die at rate $a$. Free virus is produced by  infected cells at rate  $aN$ and die at rate $k$.
This model was extended in many directions including introduction more general nonlinear terms \cite{Huang-Ma-Takeuchi-AML-2011,Korobeinikov-BMB-2007}, time delays \cite{Gourley-Kuang-Nagy-JBD-2008,Wang-Pang-Kuniya-Enatsu_AMC-2014,Wang-Liu-MMAS-2013,Rezounenko-DCDS-B-2017,Rezounenko-EJQTDE-2016,Rezounenko-DCDS-B-2018,Rezounenko-Equadiff-2017}, additional equations describing immune responses \cite{Wang-Liu-MMAS-2013,Wang-Pang-Kuniya-Enatsu_AMC-2014,Wodarz-JGV-2003,Wodarz-2007-book,Yan-Wang-DCDS-B-2012,Yousfi-Hattaf-Tridane-JMB-2011,Zhao-Xu_EJDE-2014,Rezounenko-DCDS-B-2017,Rezounenko-EJQTDE-2016} and inhomogenuous in space terms (which lead to partial differential equations models) \cite{Wang-Wang-MB-2007_HBV_spatial dependence,Wang-Huang-Zou-AA-2014,Wang-Yang-Kuniya-JMAA-2016,Rezounenko-DCDS-B-2018,Rezounenko-Equadiff-2017} (see also references therein for more information).

In this paper we, in contrast, consider a situation when viral pathogenicity is low and the infected cells are able to reproduce. Such a situation arise, for instance, when retrovirus are used for the therapy. This assumption leads to a principally different mathematical model and principally different outcomes. 
%\cite{} 
We propose the following virus dynamics model 
\begin{equation}\label{AR-2018-12-02}
\left\{
\begin{array}{l}
\frac{dC(t)}{dt}= a C(t) \left(1-b_{11}C(t) - b_{12}I(t)
\right)-\alpha\, C(t)V(t), \\
\frac{dI(t)}{dt}=a_I I(t) \left(1-b_{21}C(t) - b_{22}I(t)
\right)+\alpha\, C(t)V(t) - m I(t), \\
\frac{dV(t)}{dt}=k m I(t)-\sigma V(t), \\
\end{array}
\right.
\end{equation}
where unknowns $C(t), I(t), V(t)$ are as in the model (\ref{AR-2018-12-01})  above.

We are interested in the stability properties of an equilibrium (a stationary solution) of the model (\ref{AR-2018-12-02}). 
For the general Lyapunov stability theory see  the original work \cite{Lyapunov-1892}.

\section{Basic properties}

First, one sees that for any non-negative initial data
$$
C(0)=C_0 \ge 0, \quad I(0)=I_0 \ge 0, \quad V(0)=V_0 \ge 0,
$$
the system (\ref{AR-2018-12-02}) has a unique global (defined for all $t\ge 0$) solution. Each coordinate is non-negative for all $t\ge 0$, which is a biologically important property of the model. It follows from the standard property $\frac{dC(t)}{dt}|_{t=\tau}\ge 0$ provided $C(\tau)=0.$ Similar properties are valid for $I(t)$ and $V(t)$. In the similar standard way one shows that any solution is bounded. Moreover there is a bounded invariant region in $\mathbb{R}^3$. 

The next step of our study, is to look for possible equilibria  of the model (\ref{AR-2018-12-02}).

\subsection{Sationary solutions}

We are interested in stationary solutions of the system (\ref{AR-2018-12-02}). As a constant in-time solutions, they  satisfy 
\begin{equation}\label{AR-2018-12-03}
\left\{
\begin{array}{l}
0= a C(t) \left(1-b_{11}C(t) - b_{12}I(t)
\right)-\alpha\, C(t)V(t), \\
0=a_I I(t) \left(1-b_{21}C(t) - b_{22}I(t)
\right)+\alpha\, C(t)V(t) - m I(t), \\
0=k m I(t)-\sigma V(t). \\
\end{array}
\right.
\end{equation}
%where $...$ 

It is easy to see that there is a unique stationary solution of (\ref{AR-2018-12-03}) such that all the  coordinates are positive (inner equilibrium). In this note we are interested in this inner equilibrium and do not discuss boundary stationary solutions (when at least one coordinate is zero). 
Let us denote this unique inner  solution of (\ref{AR-2018-12-03}) by $(\widehat C, \widehat I, \widehat V)$.
%The last equation in (\ref{AR-2018-12-03}) gives
%$\widehat V = km \sigma^{-1} \widehat I$.
%We substitute it into the first and second equations of (\ref{AR-2018-12-03}) to get ...

\medskip

\begin{remark} 
We believe  this stationary solution is the most important one from the biological point of view. 
The equilibrium represents a chronic disease steady state. 
%\marginpar{???}%%%%%%%%%%%%%%%%%%%%%%%%%%%%%%%
\end{remark}

\subsection{Lyapunov stability}

Assume the parameters of the system (\ref{AR-2018-12-02}) satisfy
$$ \frac{1}{\widehat I} 
\left[ a_I b_{22} - \frac{1}{\widehat I} a_I 
\left( 1-b_{21} \widehat C - b_{22}\widehat I \right) +m 
\right] 
\cdot \frac{1}{\widehat C}
\left[ a b_{11} - \frac{1}{\widehat C}\, a \left(1-b_{11}\widehat C - b_{12}\widehat I
\right) 
%%+ \frac{1}{C} \alpha \widehat V 
\right] 
>
$$ 
\begin{equation}\label{AR-2018-12-04}
> 
\frac{1}{4} 
\left(  \frac{a b_{12}}{\widehat C}
+ \frac{b_{21}}{\widehat I}- \widehat V^2
\right)^2.
\end{equation}

Now we formulate the main result. 

\begin{theorem}\label{AR-2018-12-th1}
Let the condition (\ref{AR-2018-12-04}) be satisfied. Then there exist $\alpha_0>0, k_0>0 $ such that for any $\alpha \in (0,\alpha_0),$ and $k\in (0,k_0)$ the inner equilibrium $(\widehat C,\widehat I, \widehat V)$ of the system (\ref{AR-2018-12-02}) is locally asymptotically  stable.  
\end{theorem}

{\it Proof of theorem \ref{AR-2018-12-th1}.} 
We use the Volterra function $v(s)\equiv s - \ln\, s -1$ to construct the following Lyapunov functional

\begin{equation}\label{AR-2018-12-05}
W(C,I,V) \equiv A\cdot v\left( \frac{C}{\widehat C}\right) + 
B \cdot v\left( \frac{I}{\widehat I}\right) 
+ D \cdot v\left( \frac{V}{\widehat V}\right), 
\end{equation}
where $A,B, D$ are positive constants to be chosen below. 

It is easy to check that $W(C,I,V)>0$ for all $(C,I,V)\neq (\widehat C,\widehat I,\widehat V)$ and $W(\widehat C,\widehat I,\widehat V)=0$. 

We denote by $\dot W_{(2)}$ the derivative of $W$ along a solution of the system  (\ref{AR-2018-12-02}), which is $\dot W_{(2)}(t) = \frac{d}{dt} W(C(t),I(t),V(t))$. As usual, it is computed using the right-hand side of the system (\ref{AR-2018-12-02}) and the property $\frac{d}{dt} v(s) = 1 - s^{-1}.$

We have 
$$\dot W_{(2)} (t) 
= A \left( 1 - \frac{\widehat C}{C(t)} \right) \frac{1}{\widehat C} \left( a C(t) \left(1-b_{11}C(t) - b_{12}I(t)
\right)-\alpha\, C(t)V(t)\right)
$$
$$ +  B \left( 1 - \frac{\widehat I}{I(t)} \right) \frac{1}{\widehat I} \left( a_I I(t) \left(1-b_{21}C(t) - b_{22}I(t)
\right)+\alpha\, C(t)V(t) - m I(t)\right)
$$
\begin{equation}\label{AR-2018-12-06}
+  D \left( 1 - \frac{\widehat V}{V(t)} \right) \frac{1}{\widehat V} \left( k m I(t)-\sigma V(t)\right). 
\end{equation}
We split the above sum (\ref{AR-2018-12-06}) on there parts $\dot W_{(2)} (t) = S_1(t) + S_2(t) + S_3(t)$ to estimate them separately. 
We omit the time argument for short. 

Using $-a \widehat C \left(1-b_{11}\widehat C - b_{12}\widehat I
\right)+\alpha\, \widehat C \widehat V = 0 $, one can check that 
$$S_1 =  A  \frac{C-\widehat C}{C \widehat C} 
\left( a C \left(1-b_{11} C - b_{12} I
\right)-\alpha\,  C V
-a \widehat C \left(1-b_{11}\widehat C - b_{12}\widehat I
\right)+\alpha\, \widehat C \widehat V
\right). 
$$
Some calculations give
$$
S_1= - (C-\widehat C)^2 \cdot \frac{A}{\widehat C}
\left[ a b_{11} - \frac{1}{C}\, a \left(1-b_{11}\widehat C - b_{12}\widehat I
\right) + \frac{1}{C} \alpha \widehat V \right] 
$$
\begin{equation}\label{AR-2018-12-07}
- 
(C-\widehat C)(I-\widehat I)\cdot %\frac{1}{C}
 \frac{A a b_{12}}{\widehat C}
 - 
(C-\widehat C)(V-\widehat V)\cdot  \frac{A\alpha}{\widehat C}. 
\end{equation}
In a similar way, using $a_I \widehat I \left(1-b_{21} \widehat C - b_{22}\widehat I
\right)+\alpha\, \widehat C \widehat V - m \widehat I = 0$ (see (\ref{AR-2018-12-03})), calculations give 
$$S_2 = 
B \frac{I - \widehat I}{I \, \cdot\widehat I}
\left[
a_I I (1 - b_{21} C - b_{22} I) + \alpha C V  - m I 
\underbrace{- a_I  \widehat I (1 - b_{21}  \widehat C - b_{22}  \widehat I) - \alpha  \widehat C  \widehat V + m  \widehat I 
}%%% 
\right]
$$
$$ = B \frac{I - \widehat I}{I \, \cdot\widehat I}
\left[
a_I I \left(b_{21} (\widehat C - C) +  b_{22}(\widehat I- I)\right) + 
(I- \widehat  I) \, 
 a_I  (1 - b_{21}  \widehat C - b_{22}  \widehat I) - m  (I- \widehat  I)
 \right.
 $$
 $$\left. 
  - (\widehat V- V) C \alpha + (C- \widehat C)  \alpha \widehat V 
\right]
$$
$$
= - (I-\widehat I)^2 \cdot \frac{B}{\widehat I} 
\left[ a_I b_{22} - \frac{1}{I} a_I 
\left( 1-b_{21} \widehat C - b_{22}\widehat I \right) +m 
\right] 
$$
\begin{equation}\label{AR-2018-12-08}
= - 
(C-\widehat C)(I-\widehat I)\cdot 
\left( \frac{B\, b_{21}}{\widehat I}  - \alpha \widehat V 
\right) 
+ 
(V-\widehat V) (I-\widehat I)\cdot  \frac{C}{I} 
\frac{B \alpha}{\widehat I}. 
\end{equation}

The third step of calculations give (remind that $ \sigma  \widehat V = km {\widehat I}$ )
$$S_3=
 - 
D \frac{ V-\widehat V}{V\widehat V}\cdot 
\left[ km I - \sigma V  \underbrace{- km {\widehat I}  + \sigma  \widehat V}
\right] 
$$
\begin{equation}\label{AR-2018-12-09}
= D\, \frac{ V-\widehat V}{V\widehat V}\cdot 
\left[ km\, (I-\widehat I) - \sigma (V-\widehat V)
\right]
 = 
- (V-\widehat V)^2 \cdot  \frac{1}{V} 
\frac{D \sigma}{\widehat V} 
+ 
(V-\widehat V) (I-\widehat I)\cdot  \frac{1}{V} 
\frac{Dkm}{\widehat V}. 
\end{equation}

Hence, combining (\ref{AR-2018-12-07}), (\ref{AR-2018-12-08}), (\ref{AR-2018-12-09}),    we arrive to the form of $\dot W_{(2)} (t) = S_1(t) + S_2(t) + S_3(t)$. 
Now our goal is to find sufficient {\it local} conditions for $\dot W_{(2)} (t) < 0$ for all points (except the inner equilibrium $(\widehat C,\widehat I, \widehat V),$ where $\dot W_{(2)} (t)\equiv 0$). 
Our main idea is to compare the expression for $\dot W_{(2)}$ with auxiliary quadratic forms. Hence, knowing the conditions for the quadratic form to be positive (negative) defined, we find sufficient {\it local} conditions for $\dot W_{(2)}$ to be negative defined as a functional in a neighbourhood of $(\widehat C,\widehat I, \widehat V)$. 

Formulas (\ref{AR-2018-12-07}), (\ref{AR-2018-12-08}), (\ref{AR-2018-12-09}) are prepared to write $\dot W_{(2)}$ in the following form (notice the sign)
$$ - \dot W_{(2)} = (V-\widehat V)^2 \cdot  \omega_{11} 
 +  (I-\widehat I)^2 \cdot  \omega_{22} 
 + (C-\widehat C)^2 \cdot  \omega_{33} 
$$
\begin{equation}\label{AR-2018-12-10}
+ 
(V-\widehat V) (I-\widehat I)\cdot 2\,  \omega_{12}
+ 
(V-\widehat V) (C-\widehat C)\cdot 2\,  \omega_{13}
+ 
(C-\widehat C) (I-\widehat I)\cdot 2\,  \omega_{21}, 
\end{equation}
where all functions $\omega_{ij}$ depend on  coordinates $(C,I,V)$, parameters of the system, coordinates of the equilibrium $(\widehat C,\widehat I, \widehat V)$ and coefficients $A,B,D$. Because of the dependence of $\omega_{ij}$ on  coordinates $(C,I,V)$ the expression (\ref{AR-2018-12-10}) is not a quadratic form. We proceed as follows. Let us choose an arbitrary point $(\overline{C},\overline{I},\overline{V})$ in a small neighbourhood of $(\widehat C,\widehat I, \widehat V)$. We denote the values of $\omega_{ij}$ at this fixed point as $\overline{\omega_{ij}} \equiv \omega_{ij}(\overline{C},\overline{I},\overline{V})$ and arrive to the quadratic form
$$ \overline{W_{(2)}} = (V-\widehat V)^2 \cdot  \overline{\omega_{11}}
 +  (I-\widehat I)^2 \cdot  \overline{\omega_{22}}
 + (C-\widehat C)^2 \cdot  \overline{\omega_{33}}
$$
\begin{equation}\label{AR-2018-12-11}
+ 
(V-\widehat V) (I-\widehat I)\cdot 2\,  \overline{\omega_{12}}
+ 
(V-\widehat V) (C-\widehat C)\cdot 2\,  \overline{\omega_{13}}
+ 
(C-\widehat C) (I-\widehat I)\cdot 2\,  \overline{\omega_{23}}. 
\end{equation}
This auxiliary quadratic form $ \overline{W_{(2)}}$ is designed in such a way that $ \overline{W_{(2)}}(\overline{C},\overline{I},\overline{V}) = - \dot W_{(2)} (\overline{C},\overline{I},\overline{V})$. 
Hence, a condition to be positive defined for  $ \overline{W_{(2)}}$ implies the desired property $\dot W_{(2)} (\overline{C},\overline{I},\overline{V}) <0$. 
Since point $(\overline{C},\overline{I},\overline{V})$ is arbitrary chosen, we arrive to the local stability result.  

We apply the classical Sylvester's criterion to the quadratic form $ \overline{W_{(2)}}$. 

Let us outline the main technical idea to show that the form is positive defined about $(\widehat C,\widehat I, \widehat V)$. The  first leading principal minor of the quadratic form (\ref{AR-2018-12-10}) reads 
$\Delta_1\equiv \overline{\omega_{11}} = \frac{1}{\overline{V}} 
\frac{D \sigma}{\widehat V}>0. 
$ 
Let us consider the second leading principal minor $\Delta_2 = \overline{\omega_{11}}\, \cdot \overline{\omega_{22}} - \overline{\omega_{12}}^2. 
$
It reads 
$$\Delta_2 = \frac{1}{\overline{V}} 
\frac{D \sigma}{\widehat V} \cdot 
\frac{B}{\widehat I} 
\left[ a_I b_{22} - \frac{1}{\overline{I}} a_I 
\left( 1-b_{21} \widehat C - b_{22}\widehat I \right) +m 
\right] 
- \frac{1}{4} 
\left(
\frac{1}{\overline{V}} \cdot\frac{D km}{\widehat V}
+ \frac{\overline{C}}{\overline{I}} \cdot\frac{B \alpha}{\widehat I}
\right)^2. 
$$
We remind that we consider the {\it inner equilibrium}, so $\Delta_2\equiv \Delta_2 (\overline{C},\overline{I},\overline{V})$ is continuous in a small enough neighbourhood of $(\widehat C,\widehat I, \widehat V)$. Hence, 
the property $\Delta_2 (\widehat C,\widehat I, \widehat V)>0$ will guarantee the property 
$\Delta_2 (\overline{C},\overline{I},\overline{V})>0$ in a small enough neighbourhood. We write down the expression of  
$$
\Delta_2 (\widehat C,\widehat I, \widehat V)
= 
\frac{1}{\widehat V} 
\frac{D \sigma}{\widehat V} \cdot 
\frac{B}{\widehat I} 
\left[ a_I b_{22} - \frac{1}{\widehat I} a_I 
\left( 1-b_{21} \widehat C - b_{22}\widehat I \right) +m 
\right] 
$$
$$ 
- \frac{1}{4} 
\left(
\frac{1}{\widehat V} \cdot\frac{D km}{\widehat V}
+ \frac{\widehat C}{\widehat I} \cdot\frac{B \alpha}{\widehat I}
\right)^2 
$$
and see that %%the goal 
the property 
$\Delta_2 (\widehat C,\widehat I, \widehat V)>0$ can be reached by choosing small enough $k$ and $\alpha$. 

The same line of arguments is used to show that the third leading principal minor $\Delta_3$ is positive. We write it in a short form as follows $$\Delta_3 = \overline{\omega_{11}}\, \cdot 
( \overline{\omega_{22}} \cdot  \overline{\omega_{33}}- \overline{\omega_{23}}^2 )
 - \overline{\omega_{12}}\, \cdot
 (\overline{\omega_{12}}\cdot \overline{\omega_{33}} - \overline{\omega_{23}}\cdot \overline{\omega_{13}})
 + 
 \overline{\omega_{13}}\, \cdot
 (\overline{\omega_{12}} \cdot \overline{\omega_{23}} - \overline{\omega_{22}} \cdot \overline{\omega_{13}} ). 
$$
The detailed analysis of these three terms of $\Delta_3$ shows that the first term $\overline{\omega_{11}}\, \cdot 
( \overline{\omega_{22}} \cdot  \overline{\omega_{33}}- \overline{\omega_{23}}^2 )$ can be made positive provided condition (\ref{AR-2018-12-04}) is satisfied. Moreover, $\Delta_3$ is positive provided  $k$ and $\alpha$ are small enough. We omit the calculation here.

It completes the proof of theorem \ref{AR-2018-12-th1}.

\begin{remark} 
We notice that the smallness of parameters $k$ and $\alpha$ is not the only possible way to reach the stability. We choose this case just because of the clear biological meaning.  
\end{remark}

%\begin{remark} 
%Remark remark . ... ... \marginpar{???}%%%%%%%%%%%%%%%%%%%%%%%%%%%%%%%%%
%\end{remark}

\bigskip

{\bf Acknowledgement.} This work was supported in part by GA CR under project 16-06678S.
%%%%%%%
%%%%%%%%%%%%%%%%%%%%%%%%%%%%%%%%%%%%%%%%%%%%%%%%%%%%%%%%%%%%%%%%%%%%%%%%%%%%%%%%%%%%%%

\bigskip

\hfill December 29, 2018

\end{document}